\documentclass[10pt]{article}
\title{ An example of algebraization of analysis and Fibonacci cobweb poset characterization}
\author{Ewa Krot-Sieniawska \\
\\Institute of Computer Science,  University in Bia{\l}ystok\\
PL-15-887 Bia{\l}ystok, ul.Sosnowa 64, POLAND\\
e-mail: ewakrot@wp.pl}
\date{}

\usepackage{amsmath,amsthm}
\usepackage{amsfonts}
\chardef\bslash=`\\ 
\newtheorem{ex}{Example}[section]
\newtheorem{defn}{Definition}[section]

\newtheorem{rem}{Remark}[section]
\newtheorem{thm}{Theorem}[section]
\newtheorem{prop}{Proposition}[section]
\newtheorem{cor}{Corollary}[section]
\newtheorem{lemma}{Lemma}[section]

\begin{document}
\maketitle
\begin{abstract}
In \cite{1b,1c} inspired by O. V. Viskov \cite{23} it was shown
that the $\psi$-calculus in parts appears to be almost automatic,
natural extension of classical operator calculus of Rota - Mullin
or equivalently - of umbral calculus of Roman and Rota. At the
same time this calculus is an example of the algebraization of the
analysis - here restricted to the algebra of polynomials. The
first part of the article is the review of the recent author's
contribution \cite{2.1}. The main definitions and theorems of
Finite Fibonomial
 Operator Calculus which is a special case of $\psi$-extented Rota's finite
 operator calculus \cite{1a,1b} are   presented there. In the second part the characterization of Fibonacci Cobweb poset $P$ as DAG and oDAG
is given. The {\em dim} 2 poset such that its Hasse diagram
coincide with digraf of $P$ is constructed.
\end{abstract}


AMS 2000  numbers: 05A40, 05C20, 06A11, 11B39, 11C08

\section{Introduction}

In \cite{1d} it was shown that: ,,{\em  any $\psi $-representation
of finite operator calculus or equivalently - any $\psi
$-representation of GHW algebra makes up an example of the
algebraization of the analysis - naturally when constrained to the
algebra of polynomials.(...) Therefore the distinction in-between
difference and differentiation operators disappears. All linear
operators on $P$( the algebra  of polynomials) are both difference
and differentiation operators if the degree of differentiation or
difference operator is unlimited. For example $\frac{{d}}{{dx}} =
\sum\limits_{k \ge 1} {\frac{{d_{k}} }{{k!}}\Delta ^{k}} $ where
$d_{k} = \left[ {\frac{{d}}{{dx}}x^{\underline {k}} } \right]_{x =
0}  = \left( { - 1} \right)^{k - 1}\left( {k - 1} \right)!$ or
$\Delta =\sum\limits_{n \ge 1} {\frac{{\delta _{n}} }{{n!}}}
\frac{{d^{n}}}{{dx^{n}}}$ where $\delta _{n} = \left[ {\Delta
x^{n}} \right]_{x = 0} =1$. Thus the difference and differential
operators and equations are treated on the same footing.}'' The
authors  goal  was there ,,{\em to  deliver the general scheme of
"$\psi $-umbral" algebraization of the analysis of general
differential operators \cite{12}.}'' One may consider a plenty of
different special cases of these $\psi$-extensions. Each of them
can be obtained by the special choice of an admissible sequence
$\psi = \left\{ {\psi _{n} \left( {q} \right)} \right\}_{n \ge 0}
$ ; $\psi _{n} \left( {q} \right) \ne 0$; $n \ge 0$ and $\psi _{0}
\left( {q} \right) = 1$, or equivalently, by the special choice of
an sequence $n_{\psi}  $, where (\cite{1, 1a, 1b, 1d})
$$
n_{\psi}  \equiv \psi _{n - 1} \left( {q} \right)\psi _{n}^{-1}
\left( {q} \right),n \geq 0,$$ and each of them constitutes  a
representation of GHW algebra and  provides  an example of the
algebraization of the analysis.

One of most interesting cases is the  so called Finite Fibonomial
Operator Calcus (FFOC). Its idea comes from \cite{1a} and it was
considered by the present author in \cite{2.1}. It is also the
main object of this work. In the first part of it, we present some
definitions and theorems of FFOC. We take there
$n_{\psi}=n_F=F_n$, where the famous Fibonacci
sequence$\{F_{n}\}_{n \geq 0}$
$$\Bigg\{\begin{array}{l}F_{n+2}=F_{n+1}+F_{n}\\
F_{0}=0,\;F_{1}=1\end{array}$$ is attributed and refered to the
first edition (lost) of "Liber Abaci" (1202) by Leonardo Fibonacci
(Pisano)(see edition from 1228 reproduced as "Il Liber Abaci di
Leonardo Pisano publicato secondo la lezione Codice Maglibeciano
by Baldassarre Boncompagni in Scritti di Leonardo Pisano" , vol.
1, (1857) Rome).

In order to formulate main results of FFOC the following objects
for the sequence $F=\left\{F_{n}\right\}_{n \geq 0}$ are defined:
\renewcommand{\labelenumi}{(\arabic{enumi})}
\begin{enumerate}
\item $F$-factorial:
$$F_{n}!=F_{n}F_{n-1}...F_{2}F_{1},\;\;\;F_{0}!=1.$$
\item $F$-binomial (Fibonomial ) coefficients \cite{4}:
$$\binom{n}{k}_{F}=\frac{n_{F}^{\underline{k}}}{k_{F}!}=
\frac{F_{n}F_{n-1}\ldots F_{n-k+1}}{F_{k}F_{k-1}\ldots
F_{2}F_{1}}=
\frac{F_{n}!}{F_{k}!F_{n-k}!},\;\;\;\binom{n}{0}_{F}=1.$$
\end{enumerate}
It is known that $\binom{n}{k}_{F}\in {\bf N}$ for every $n,k \in
{\bf N}\cup {0}$.

In \cite{13} some applications of $\psi$ -extensions of the umbral
calculus (including FFOC) were presented. As annouced there,  the
combinatorial interpretation of fibonomial coefficients has been
found by A. K. Kwa\'sniewski in \cite{4.1}. It was done by the use
of the so called Fibonacci cobweb poset \cite{3.1, 10.1, 11.1,
4.1, 1c, 55f}. In \cite{2.3, 2.4} the incidence algebra of the
Fibonacci cobweb poset was considered by the present author. In
the second part of this work the characterisation of this poset
$P$ as DAG and oDAG is given. The {\em dim} 2 poset such that its
Hasse diagram coincide with digraf of $P$ is constructed. Directed
acyclic graphs  (DAGs) have many important applications in
computer science, including: the parse tree constructed by a
compiler, a reference graph that can be garbage collected using
simple reference counting, dependency graphs such as those used in
instruction scheduling and makefiles, dependency graphs between
classes formed by inheritance relationships in object-oriented
programming languages. In theoretical physics a directed acyclic
graph can be used to represent spacetime as a causal set. In
bioinformatics, DAGs can be used to find areas of synteny between
two genomes. They can be also used  in abstract process
descriptions such as workflows and some models of provenance.

\section{Finite Fibonomial Operator Calculus}

\subsection{Operators and polynomial sequences} Let {\bf P} be the
algebra of polynomials over the field ${\mathbb{K}}$ of
characteristic zero.
\begin{defn}
 The linear operator   $\partial_{F}:{\bf P}\rightarrow{\bf P}$   such that
   $\partial_{F}x^{n}=F_{n}x^{n-1}$  for $n \geq 0$ is named the
   $F$-derivative.
\end{defn}
\begin{defn}
The $F$-translation operator is the linear operator \\
$E^{y}(\partial_{F}): {\bf P}\rightarrow{\bf P}$ of the form:
$$E^{y}(\partial_{F})=\exp_{F}\{ y \partial_{F} \}=
\sum_{k \geq 0} \frac{y^{k} \partial_{F}^{k}}{F_{k}!},\;\;\;\;\;\;
y\in{\bf K}$$
\end{defn}
\begin{defn}
$$\forall_{p\in{\bf P}}\;\;\;\;p(x+_{F}y)=E^{y}(\partial_{F})p(x)\;\;\;\;
x,y\in{\bf K}$$
\end{defn}
 \begin{defn}
A linear operator $T:{\bf P}\rightarrow{\bf P}$ is said to be
 $\partial_{F}$-shift invariant iff
$$\forall_{y \in {\bf K}}\;\;\;\; [T,E^{y}(\partial_{F})]=
TE^{y}(\partial_{F})-E^{y}(\partial_{F})T=0$$ We shall denote by
$\Sigma_{F}$ the algebra of $F$-linear $\partial_{F}$-shift
invariant operators.
\end{defn}

 \begin{defn}
Let $Q(\partial_{F})$ be a formal series in powers of
$\partial_{F}$ and $Q(\partial_{F}):{\bf P}\rightarrow{\bf P}$.
$Q(\partial_{F})$ is said to be $\partial_{F}$-delta operator iff
\renewcommand{\labelenumi}{(\alph{enumi})}
\begin{enumerate}
\item $Q(\partial_{F}) \in \Sigma_{F}$ \item $Q(\partial_{F})(x) =
const \neq 0$
 \end{enumerate}
 \end{defn}

 Under quite natural specification the proofs of most statements
 might be reffered to \cite{1}(see also references therein).

 The particularities of the case considered here are revealed in
 the sequel. There the scope of
 new possibilities is initiated by means of unknown before
 examples.
\begin{prop}
Let $Q(\partial_{F})$ be the $\partial_{F}$-delta operator. Then
$$\forall_{ c\in{\bf K}}\;\;\;\;Q(\partial_{F})c=0.$$
\end{prop}
\begin{prop}
Every $\partial_{F}$-delta operator reduces degree of any
polynomial by one.
\end{prop}

\begin{defn}
The polynomial sequence $\{q_{n}(x)\}_{n\geq 0}$ such that \\
$deg\;q_{n}(x)=n$ and:
\renewcommand{\labelenumi}{(\arabic{enumi})}
\begin{enumerate}
\item $q_{0}(x)=1;$ \item $q_{n}(0)=0,\; n\geq 1;$ \item
$Q(\partial_{F})q_{n}(x)=F_{n}q_{n-1}(x),\;\;n \geq 0$
\end{enumerate}
 is called $\partial_{F}$-basic polynomial sequence of
the $\partial_{F}$-delta operator $Q(\partial_{F})$.
\end{defn}
\begin{prop}
For every $\partial_{F}$-delta operator $Q(\partial_{F})$ there
exists the uniquely determined $\partial_{F}$-basic polynomial
sequence $\{q_{n}(x)\}_{n\geq 0}$.
\end{prop}
\begin{defn}
 A polynomial sequence $\{p_{n}(x)\}_{n\geq 0}$ ($deg\;p_{n}(x)=n$)
 is of $F$-binomial (fibonomial) type if it satisfies the condition
 $$E^{y}(\partial_{F})p_{n}(x)=p_{n}(x+_{F}y)=\sum_{k\geq 0} \binom{n}{k}_{F}
 p_{k}(x)p_{n-k}(y)\;\;\;\forall_{  y\in {\bf K}}$$
\end{defn}
\begin{thm}
The  polynomial sequence $\{p_{n}(x)\}_{n\geq 0}$ is a
$\partial_{F}$-basic
 polynomial sequence of some $\partial_{F}$-delta operator $Q(\partial_{F})$
 iff it is a sequence of\\ $F$-binomial type.
\end{thm}
\begin{thm} {\em (First Expansion Theorem)}\\
Let $T \in \Sigma_{F}$ and let $Q(\partial_{F})$ be a
$\partial_{F}$-delta operator with $\partial_{F}$-basic polynomial
sequence $\{q_{n}\}_{n \geq 0}$. Then
$$T = \sum _{n \geq 0} \frac{a_{n}}{F_{n}!}Q(\partial_{F})^{n};
\quad a_{n} = [Tq_{k}(x)]_{x=0}.$$
\end{thm}

\begin{thm}{\em (Isomorphism Theorem)}\\
Let $\Phi_{F}={\bf K}_{F}[[t]]$ be the algebra of formal {\em
exp}$_{F}$ series in $t \in {\bf K}$ ,i.e.:
$$f_{F}(t)\in \Phi_{F}\;\;\;\;  iff\;\;\;\; f_{F}(t)= \sum_{k \geq 0}
\frac{a_{k}t^{k}}{F_{k}!}\;\;\; for\;\;\; a_{k}\in {\bf K},$$ and
let the $Q(\partial_{F})$ be a $\partial_{F}$-delta operator. Then
$\Sigma_{F} \approx \Phi_{F}$.
 The isomorphism \\ $\phi : \Phi_{F} \rightarrow \Sigma_{F}$ is given by the
 natural correspondence:
 $$f_{F}  \left( {t} \right) = \sum_{k \geq 0}
\frac{a_{k}t^{k}}{F_{k}!}\; \buildrel {into} \over \longrightarrow
\; T_{\partial_{F}} = \sum_{k \geq 0}
\frac{a_{k}}{F_{k}!}Q(\partial_{F})^{k}. $$
\end{thm}
\begin{rem} {\em
In the algebra $\Phi_{F}$ the product is given by the fibonomial
convolution, i.e.:
$$
\left( {\;\sum_{k \geq 0} {\frac{{a_{k}} }{{F_{k}  !}}} x^{k}\;}
\right) \left( {\;\sum_{k \geq 0} {\frac{{b_{k}} }{{F_{k}  !}}}
x^{k}\;} \right)= \left( {\;\sum_{k \geq 0} {\frac{{c_{k}}
}{{F_{k}  !}}} x^{k}\;} \right)$$ where
$$c_{k} = \sum_{l \geq 0} \binom{k}{l}_{F}  a_{l} b_{k -l}.$$ }
\end{rem}
\begin{cor}
Operator $T \in \Sigma_{F}$ has its inverse $T^{ - 1}\; \in
\Sigma_{\psi}$   iff    $T1 \neq 0$.
\end{cor}
\begin{rem}{\em
 The $F$-translation operator
$E^{y}\left( {\partial _{F} }  \right) = \exp_{F}  \{ y\partial
_{F} \} $ is invertible in $\Sigma_{F}$ but it is not a $\partial
_{F} $-delta operator.
 No one of $\partial_{F}  $-delta operators $Q\left( {\partial _{F} }  \right)$
is invertible with respect to the formal series "F-product".  }
\end{rem}
\begin{cor}
Operator $R(\partial _{F} ) \in \Sigma_{F} $ is a $\partial
_{F}$-delta operator iff $a_{0} = 0$ and $a_{1} \neq 0$, where
$R(\partial _{F} )  = \sum_{n \geq 0} \frac{{a_{n}} }{F_{n} !}
Q\left( \partial _{F}  \right)^{n}$ or equivalently : $r(0) = 0$
{\em \&} $r'(0) \neq 0$ where $r(x) = \sum\limits_{k \geq
0}\frac{a_{k}}{F_{k}!}x^{k}\;$ is the correspondent of $R(\partial
_{F})\;$ under the Iomorphism Theorem.
\end{cor}
\begin{cor}
Every $\partial_{F}  $-delta operator $Q\left( \partial _{F}
\right)$ is a function $Q(\partial _{F})  $ according to the
expansion
$$Q\left( \partial _{F}   \right) = \sum\limits_{n \geq 1}\frac{q_{n}
}{F_{n}  !} \partial _{F} ^{n}$$ This $F$-series will be called
the $F$-indicator of the $Q(\partial_{F})$.
\end{cor}
\begin{rem}{\em
$\exp_{F} \{zx\}$ is the $F$-exponential generating function for
\\ $\partial _{F}  $-basic polynomial sequence $\left\{ x^{n}
\right\}_{n = 0}^{\infty}  $ of the $\partial _{F}  $ operator.}
\end{rem}

\begin{cor}
The $F $-exponential generating function for $\partial _{F}
$-basic polynomial sequence $\left\{ p_{n} \left( {x} \right)
\right\}_{n = 0}^{\infty}  $ of the $\partial _{F}  $-delta
operator $Q\left( {\partial _{F} }  \right)$ is given by the
following formula
$$\sum\limits_{k \geq 0} \frac{{p_{k} \left( x \right)}}{F_{k}  !}
z^{k}\; = \exp_{F}  \{ xQ^{ - 1}\left( z \right)\}$$ where
$$Q \circ Q^{-1}=Q^{-1} \circ Q=I=id.$$
\end{cor}
 \begin{ex} \label{deltaop}

  {\em The following operators are the examples of $\partial_{F}$-delta\\
  operators:}
\renewcommand{\labelenumi}{\em (\arabic{enumi})}
\begin{enumerate} {\em
\item $\partial_{F}$; \item $F$-difference operator
$\Delta_{F}=E^{1}(\partial_{F})-I$
 such that\\
$(\Delta_{F}p)(x)=p(x+_{F}1)-p(x)$ for every $p \in {\bf P}$ ;
\item The operator $\nabla_{F}=I-E^{-1}(\partial_{F})$ defined as follows:\\
$(\nabla_{F}p)(x)=p(x)-p(x-_{F}1)$ for every $p \in {\bf P}$;
\item $F$-Abel operator:
$A(\partial_{F})=\partial_{F}E^{a}(\partial_{F})= \sum\limits_{k
\geq 0} \frac{a^{k}}{F_{k}!}\partial_{F}^{k+1}$; \item
$F$-Laguerre operator of the form:
$L(\partial_{F})=\frac{\partial_{F}}
{\partial_{F}-I}=\sum\limits_{k \geq 0}\partial_{F}^{k+1}$.}
\end{enumerate}
\end{ex}
\begin{defn}
The $\hat{x}_{F}$-operator is the linear map $\hat{x}_{F}:{\bf P}
\rightarrow {\bf P}$ such that\\
$\hat{x}_{F}x^{n}=\frac{n+1}{F_{n+1} }x^{n+1} \; for \;\; n\geq
0$. ($\;[\partial_{F},\hat{x}_{F}]=id$.)
\end{defn}
\begin{defn}
A linear map {\bf '} : $\Sigma _{F} \rightarrow \Sigma _{F} $ such that   \\
$T\;${\bf '} = $T\;\hat {x}_{F}  - \hat {x}_{F}  T$ =
\textbf{[}$T$, $\hat
 {x}_{F} $\textbf{]}\\
 is called the Graves-Pincherle $F$-derivative {\em \cite{9,10}}.
\end{defn}
\begin{ex} \textrm{ }
 \renewcommand{\labelenumi}{(\arabic{enumi})}
\begin{enumerate}
\item $\partial_{F}${\bf '}=$I=id$; \item $(\partial_{F})^{n}${\bf
'}=$n\partial_{F}^{n-1}$
\end{enumerate}
\end{ex}

 According to the example above the Graves-Pincherle
$F$-derivative is the formal derivative with respect to
$\partial_{F}$ in $\Sigma_{F}$ i.e., $T${\bf }'$\;(\partial_{F})
\in \Sigma_{F}$ for any $T \in \Sigma_{F} $.
\begin{cor}\label{corfourone}
Let $t\left( {z} \right) $ be the indicator of operator $T \in
\Sigma _{F} $. Then\\ $t'\left( {z} \right)$ is the indicator of
$T${\bf '}$ \in \Sigma _{F}$.
\end{cor}

Due to the isomorphism theorem and the Corollaries above the
Leibnitz rule
 holds .
\begin{prop} \label{propfourone}
($T S$)\textbf{'} $=T$\textbf{'} $S\; + ST$\textbf{'}$\;$ ; $T$,
$S \in \; \Sigma _{F} $.
\end{prop}

As an immediate consequence of the Proposition~\ref{propfourone}
we get
\begin{center}
($S  ^{n}\;$)\textbf{'}= n $S$\textbf{'}$S^{n - 1}\; \quad
\forall_{ S \in \Sigma _{F}}$.
\end{center}
From the isomorphism theorem we insert that the following is true.

\begin{prop}
$Q\left( {\partial _{F} }  \right)$ is the $\partial _{F}$-delta
operator iff
 there exists invertible
$S\in \Sigma_{F}$ such that
$$Q\left( {\partial _{F} }  \right) \; = \;
\partial _{F} S.$$
\end{prop}

The Graves-Pincherle $F$-derivative notion appears very effective
while formulating expressions for $\partial _{F}  $-basic
polynomial sequences of the given $\partial _{F}  $-delta operator
$Q\left( {\partial _{F} } \right)$.
  \begin{thm} \label{lr}
{\em ($F$-Lagrange and $F$-Rodrigues formulas) \cite{1,6,12}}\\
Let $\{q_{n}\}_{n \geq 0}$ be $\partial_{F}$-basic sequence of the
delta operator $Q(\partial_{F})$, $Q(\partial_{F})=\partial_{F}P$
($P \in \Sigma_{F}$, invertible). Then for $n\geq0$:

\begin{enumerate}
\renewcommand{\labelenumi}{\em (\arabic{enumi})}
\item\label{one} $q_{n}(x) = Q\left( \partial _{F}
\right)$\textbf{'}
 $P^{-n-1}\;x^{n}$ ;

\item\label{two} $q_{n}(x) = P^{-n} x^{n} - \frac{F_{n}}{{n}}$
($P^{ - n}\;$) \textbf{'}$x^{n-1};$

\item\label{three} $q_{n}(x) = \frac{{F_{n} } }{{n}}\hat
{x}_{F}P^{ - n} x^{n-1}$;

\item\label{four} $q_{n}(x) = \frac{F_{n}}{n}\hat {x}_{F}(Q\left(
{\partial _{F} }  \right)$\textbf{'} )$^{-1} q_{n-1}(x)$  {\em
($\leftarrow$ Rodrigues $F $-formula )}.
\end{enumerate}
\end{thm}
\begin{cor}
Let $Q(\partial_{F})=\partial_{F}S$ and
$R(\partial_{F})=\partial_{F}P$ be
 the $\partial_{F}$-delta operators with the $\partial_{F}$-basic sequences
 $\{q_{n}(x)\}_{n \geq 0}$ and $\{r_{n}(x)\}_{n \geq 0}$ respectively. Then:
 \renewcommand{\labelenumi}{\em (\arabic{enumi})}
\begin{enumerate}
\item
$q_{n}(x)=R$\textbf{'}$(Q$\textbf{'}$)^{-1}S^{-n-1}P^{n+1}r_{n}(x),
\;\;\;n \geq 0$; \item
$q_{n}(x)=\hat{x}_{F}(PS^{-1})^{n}\hat{x}_{F}^{-1}r_{n}(x),\;\;\;
n>0$.
\end{enumerate}
\end{cor}
The formulas of the Theorem~\ref{lr} can be used to find
$\partial_{F}$-basic
 sequences of the $\partial_{F}$-delta operators from the
  Example~\ref{deltaop}.
\begin{ex} \textrm{  }
\renewcommand{\labelenumi}{\em (\arabic{enumi})}
\begin{enumerate} {\em
\item The polynomials $x^{n},\;n \geq 0$ are $\partial_{F}$-basic
for
 $F$-derivative $\partial_{F}$.
 \item Using Rodrigues formula in a straighford  way one can find the following first
 $\partial_{F}$-basic polynomials
 of the operator $\Delta_{F}$:\\
 $q_{0}(x)=1\\
 q_{1}(x)=x\\
 q_{2}(x)=x^{2}-x\\
 q_{3}(x)=x^{3}-4x^{2}+3x\\
 q_{4}(x)=x^{4}-9x^{3}+24x^{2}-16x\\
 q_{5}(x)=x^{5}-20x^{4}+112.5x^{3}-250x^{2}+156.5x\\
 q_{6}(x)=x^{6}-40x^{5}+480x^{4}-2160x^{3}+4324x^{2}-2605x .$
 \item Analogously to the above example we find the following first
  $\partial_{F}$-basic polynomials of the operator
 $\nabla_{F}$:\\
$q_{0}(x)=1\\
 q_{1}(x)=x\\
 q_{2}(x)=x^{2}+x\\
 q_{3}(x)=x^{3}+4x^{2}+3x\\
 q_{4}(x)=x^{4}+9x^{3}+24x^{2}+16x\\
 q_{5}(x)=x^{5}+20x^{4}+112.5x^{3}+250x^{2}+156.5x\\
 q_{6}(x)=x^{6}+40x^{5}+480x^{4}+2160x^{3}+4324x^{2}+2605x .$
 \item Using Rodrigues formula in a straighford  way one finds the following first
 $\partial_{F}$-basic polynomials of $F$-Abel operator:\\
$A^{(a)}_{0,F}(x)=1\\
A^{(a)}_{1,F}(x)=x\\
A^{(a)}_{2,F}(x)=x^{2}+ax\\
A^{(a)}_{3,F}(x)=x^{3}-4ax^{2}+2a^{2}x\\
A^{(a)}_{4,F}(x)=x^{4}-9ax^{3}+18a^{2}x^{2}-3a^{3}x .$

  \item In order to find $\partial_{F}$-basic polynomials of
   $F$-Laguerre operator $L(\partial_{F})$  we use formula (3) from Theorem~\ref{lr}:
   \begin{multline*}
   L_{n,F}(x)=\frac{F_{n}}{n}\hat{x}_{F}\left(\frac{1}{\partial_{F}-1}
   \right)^{-n}x^{n-1}=\frac{F_{n}}{n}\hat{x}_{F}(\partial_{F}-1)^{n}
   x^{n-1}=\\=\frac{F_{n}}{n}\hat{x}_{F}\sum_{k=0}^{n}(-1)^{k}\binom{n}{k}
   \partial_{F}^{n-k}x^{n-1}=\frac{F_{n}}{n}\hat{x}_{F}\sum_{k=0}^{n}
   (-1)^{k}\binom{n}{k}(n-1)^{\underline{n-k}}_{F} x^{k-1}=\\=
   \frac{F_{n}}{n}\sum_{k=1}^{n}(-1)^{k}\binom{n}{k}(n-1)^{\underline
   {n-k}}_{F} \frac{k}{F_{k}}x^{k}.
  \end{multline*}}
  \end{enumerate}
 \end{ex}

 \subsection{Sheffer $F$-polynomials}
\begin{defn}
A polynomial sequence $\{s_{n}\}_{n\geq 0}$ is called the sequence
of Sheffer $F$-polynomials of the $\partial_{F}$-delta operator
$Q(\partial_{F})$ iff

\renewcommand{\labelenumi}{\em (\arabic{enumi})}
\begin{enumerate}
\item $s_0(x)=const\neq 0$ \item
$Q(\partial_{F})s_{n}(x)=F_{n}s_{n-1}(x);\; n\geq 0.$
\end{enumerate}
\end{defn}

\begin{prop} \label{shefprop}
Let $Q(\partial_{F})$ be $\partial_{F}$-delta operator with
$\partial_{F}$-basic polynomial sequence $\{q_{n}\}_{n  \geq 0}$.
Then $\{s_{n}\}_{n\geq 0}$ is the sequence of Sheffer
$F$-polynomials of $Q(\partial_{F})$ iff there exists an
invertible $S \in \Sigma _{F}$ such that $s_{n}(x)=S^{-1}q_{n}(x)$
for  $n\geq 0$. We shall refer to a given labeled by
$\partial_{F}$-shift invariant invertible operator $S$   Sheffer
$F$-polynomial sequence $\{s_{n}\} _{n\geq 0}$ as the sequence of
Sheffer $F$-polynomials of the $\partial_{F}$-delta operator
$Q(\partial_{F})$ relative to $S$.
\end{prop}
\begin{thm}\label{thfourthree} {\em (Second $F $- Expansion Theorem)}\\
Let $Q\left( {\partial _{F} }  \right)$ be the $\partial _{F}
$-delta operator $Q\left( {\partial _{F} }  \right)$ with the
$\partial _{F} $-basic polynomial sequence $\left\{ {q_{n} \left(
{x} \right)} \right\}_{n\geq 0}  $. Let $S$ be an
\textit{invertible} $\partial _{F}  $-shift invariant operator and
let $\left\{ {s_{n} \left( {x} \right)} \right\}_{n \geq 0} $ be
its sequence of Sheffer $F $-polynomials. Let $T$ be \textit{any}
$\partial _{F}  $-shift invariant operator and let \textit{p(x)}
be any polynomial. Then the following identity holds :

\begin{center}
$\forall_{ y \in K} \wedge \; \forall_{ p \in P} \; \quad
(Tp)\left( {x + _{F} y} \right)
=\left[E^{y}(\partial_{F})p\right](x)=T \sum\limits _{k \geq
0}\frac{s_{k} \left( y \right)}{F_{k}  !} Q\left( {\partial _{F} }
\right)^{k}S\;  Tp\left( {x} \right)$ .
\end{center}
\end{thm}
\begin{cor}
Let ${s_{n}(x)}_{n \geq 0}$ be a sequence of Sheffer
$F$-polynomials of  a  $\partial_{F}$-delta operator
$Q(\partial_{F})$ relative to $S$.Then:
$$S^{-1}=\sum_{k \geq 0}\frac{s_{k}(0)}{F_{k}!}Q(\partial_{F})^{k}.$$
\end{cor}
\begin{thm}({\em The Sheffer $F$-Binomial Theorem})\\
Let $Q(\partial_{F})$, invertible $S \in\Sigma_{F},{q_{n}(x)} _{n
\geq 0},{s_{n}(x)}_{n \geq 0}$ be as above. Then:
 \begin{center}
$$E^{y}(\partial_{F})s_{n}(x)=s_{n}(x+_{F}y)=\sum_{k \geq 0}\binom{n}{k}_{F}
s_{k}(x)q_{n-k}(y).$$
\end{center}
\end{thm}
\begin{cor}
$$s_{n}(x)=\sum_{k \geq 0}\binom{n}{k}_{F}s_{k}(0)q_{n-k}(x)$$
\end{cor}
\begin{prop}
Let $Q\left( {\partial _{F} }  \right)$ be a $\partial _{F}
$-delta operator. Let $S$ be an invertible $\partial _{F } $-shift
invariant operator. Let $\left\{ {s_{n} \left( {x} \right)}
\right\}_{n \geq 0} $ be a polynomial sequence. Let
\begin{center}
$\forall_{ a \in K }\wedge \; \forall_{ p \in P} \quad E^{a}\left(
{\partial _{F} }  \right)p\left( {x} \right) = \sum\limits_{k \ge
0} {\frac{{s_{k} \left( {a} \right)}}{{F_{k}  !}}} Q\left(
{\partial _{F} }  \right)^{k}S_{\partial _{F} }  \;p\left( {x}
\right)$ .
\end{center}
Then the polynomial sequence $\left\{ {s_{n} \left( {x} \right)}
\right\}_{n \geq 0}$ is the sequence of Sheffer  $F $-polynomials
of the $\partial _{F}  $-delta operator $Q\left( {\partial _{F} }
\right)$ relative to $S$.
\end{prop}
\begin{prop}
Let $Q\left( {\partial _{F} }  \right)$and  $S$ be as above. Let
\textit{q(t)} and s\textit{(t)}  be the indicators of $Q\left(
{\partial _{F} }  \right)$ and $S$ operators. Let
\textit{q$^{-1}$(t )} be the inverse $F $-exponential formal power
series inverse to \textit{q(t)}. Then the $F $-exponential
generating function of Sheffer $F $-polynomials sequence $\left\{
{s_{n} \left( {x} \right)} \right\}_{n \geq 0}$ of $Q\left(
{\partial _{F} }  \right)$ relative to $S\;$is given by
\[
\;\sum\limits_{k \ge 0} {\frac{{s_{k} \left( {x} \right)}}{{F_{k}
!}}} z^{k}\; = \;\left(s\left( {q^{ - 1}\left( {z} \right)}
\right)\right)^{-1}\;\exp_{F}  \{ xq^{ - 1}\left( {z} \right)\}.
\]
\end{prop}
\begin{prop}
A sequence $\left\{ {s_{n} \left( {x} \right)} \right\}_{n \geq
0}$ is the sequence of Sheffer  $F $-polynomials of the $\partial
_{F} $-delta operator $Q\left( {\partial _{F} }  \right)$ with the
$\partial _{F}  $-basic polynomial sequence $\left\{ {q_{n} \left(
{x} \right)} \right\}_{n \geq 0}$ iff
\[
s_{n} \left( {x + _{F}  y} \right)= \sum\limits_{k \ge 0}
\binom{n}{k}_{F}  s_{k} \left( {x} \right)q_{n - k} \left( {y}
\right).
\]
for all $y \in {\bf K}$
\end{prop}

\begin{ex}{\em
Hermite $F$-polynomials are Sheffer $F$-polynomials of the\\
$\partial_{F}$ -delta operator $\partial_{F}$ relative to
invertible $S \in \Sigma_{F}$ of the form \\ $S=\exp_{F} \{
\frac{a \partial_{F}^{2}}{2}\}$. One can get them by formula (see
Proposition~\ref{shefprop} ):
$$H_{n,F}(x)=S^{-1}x^{n}=\sum\limits_{k \geq 0}\frac{(-a)^{k}}{2^{k}F_{k}!}
n^{\underline{2k}}_{F}x^{n-2k}.$$}
\end{ex}
\begin{ex}   {\em
 Let $S=(1-\partial_{F})^{-\alpha- 1}$. The Sheffer $F$-polynomials of\\
 $\partial_{F}$-delta operator $L(\partial_{F})=\frac{\partial_{F}}
 {\partial_{F}-1}$  relative to $S$ are Laguerre $F$-polynomials of order
 $\alpha$ . By Proposition~\ref{shefprop} we have
 $$L^{(\alpha)}_{n,F}=(1-\partial_{F})^{\alpha+1}L_{n,F}(x),$$
From the above formula and using Graves-Pincherle $F$-derivative
we get
$$L^{(\alpha)}_{n,F}(x)=\sum\limits_{k \geq 0}\frac{F_{n}!}{F_{k}!}\binom
{\alpha+n}{n-k}(-x)^{k}$$ for $\alpha \neq -1$}.
\end{ex}
\begin{ex} {\em
Bernoullie's $F$-polynomials of order 1 are Sheffer
$F$-polynomials of \\$\partial_{F}$ -delta operator $\partial_{F}$
related to invertible $S=\left(\frac
{\exp_{F}\{\partial_{F}\}-I}{\partial_{F}}\right)^{-1}$. Using\\
Proposition~\ref{shefprop} one arrives at
 \begin{multline*}
 B_{n,F}(x)=S^{-1}x^{n}=\sum_{k \geq 1}\frac{1}{F_{k}!}\partial_{F}
^{k-1}x^{n}=\sum_{k \geq
1}\frac{1}{F_{k}}\binom{n}{k-1}_{F}x^{n-k+1}=\\= \sum_{k \geq
0}\frac{1}{F_{k+1}}\binom{n}{k}_{F}x^{n-k}
\end{multline*}}
\end{ex}

\begin{thm}\label{recrel}
{\em (Reccurence relation for Sheffer $F$-polynomials)}\\
Let $Q,S,\{s_{n}\}_{n\geq 0} $ be as above. Then the following
reccurence formula holds:
$$
s_{n+1}(x)=\frac{F_{n+1}}{n+1}\left[\hat{x}_{F}-\frac{S'}{S}\right]
\left[Q(\partial_{F})'\right]^{-1}s_{n}(x);\; n\geq 0.
$$
\end{thm}
\begin{ex}
{\em The reccurence formula for the Hermite $F$-polynomials is:}
$$H_{n+1,F}(x)=\hat{x}_{F}H_{n,F}(x)-\hat{a}_{F}F_{n}H_{n-1,F}(x)$$
\end{ex}
\begin{ex}
{\em The reccurence relation for the Laguerre $F$-polynomials is:}
\begin{multline*}
L_{n+1,F}^{(\alpha)}(x)=-\frac{F_{n+1}}{n+1}[\hat{x}_{F}-(\alpha+1)(1-\partial_{F})^{-1}
](\partial_{F}-1)^{2}L_{n,F}^{(\alpha)}(x)\\=\frac{F_{n+1}}{n+1}[\hat{x}_{F}(\partial_{F}-1)
+\alpha+1]L_{n,F}^{(\alpha+1)}(x).
\end{multline*}
\end{ex}
 \subsection{Some examples of $F$-polynomials}
  \renewcommand{\labelenumi}{(\arabic{enumi})}
\begin{enumerate}
\item Here are the examples of Laguerre $F$-polynomials of order
$\alpha=-1$:
\\ \textrm{}\\
$L_{0,F}(x)=1\\ \textrm{ } \\
L_{1,F}(x)=-x\\ \textrm{ } \\
L_{2,F}(x)=x^{2}-x\\ \textrm{ } \\
L_{3,F}(x)=-x^{3}+4x^{2}-2x\\ \textrm{ } \\
L_{4,F}(x)=x^{4}-9x^{3}+18x^{2}-6x\\  \textrm{ } \\
L_{5,F}(x)=-x^{5}+20x^{4}-905x^{3}+1280x^{2}-30x\\  \textrm{ } \\
L_{6,F}(x)=x^{6}-40x^{5}+400x^{4}-1200x^{3}+1200x^{2}-
240x\\  \textrm{ } \\
L_{7,F}(x)=-x^{7}+78x^{6}-1560x^{5}+
10400x^{4}-23400x^{3}+18720x^{2}-\\
 \textrm{}\\\;\;\;\;\;\;-3120x\\ \textrm{}\\
L_{8,F}(x)=x^{8}-147x^{7}+5733x^{6}-76440x^{5}+
382200x^{4}-687960x^{3}+\\ \;\;\;\;+458640x^{2}-65520x$ \item Here
are the examples of Laguerre $F$-polynomials of order $\alpha=1$:
\\ \textrm{}\\
$L^{(1)}_{0,F}(x)=1\\ \textrm{ } \\
L^{(1)}_{1,F}(x)=-x+2\\ \textrm{ } \\
L^{(1)}_{2,F}(x)=x^{2}-3x+3\\ \textrm{ } \\
L^{(1)}_{3,F}(x)=-x^{3}+8x^{2}-12x+8\\ \textrm{ } \\
L^{(1)}_{4,F}(x)=x^{4}-15x^{3}+60x^{2}-60x+30\\  \textrm{ } \\
L^{(1)}_{5,F}(x)=-x^{5}+30x^{4}-225x^{3}+600x^{2}-450x+240\\
 \textrm{ } \\
L^{(1)}_{6,F}(x)=x^{6}-56x^{5}+840x^{4}-4200x^{3}+8400x^{2}-
5040x+1680$ \item  Here we give some examples of the Bernoullie's
$F$-polynomials of
order 1:\\
\textrm{}\\
$B_{0,F}(x)=1\\ \textrm{ } \\
B_{1,F}(x)=x+1\\ \textrm{ } \\
B_{2,F}(x)=x^{2}+x+\frac{1}{2}\\ \textrm{ } \\
B_{3,F}(x)=x^{3}+2x^{2}+x+\frac{1}{3}\\ \textrm{ } \\
B_{4,F}(x)=x^{4}+3x^{3}+3x^{2}+x+\frac{1}{5}\\  \textrm{ } \\
B_{5,F}(x)=x^{5}+5x^{4}+\frac{15}{2}x^{3}+5x^{2}+x+\frac{1}{8}\\ \textrm{ } \\
B_{6,F}(x)=x^{6}+8x^{5}+20x^{4}+20x^{3}+8x^{2}+x+\frac{1}{13}\\  \textrm{ } \\
B_{7,F}(x)=x^{7}+13x^{6}+52x^{5}+\frac{260}{3}x^{4}+52x^{3}+13x^{2}+x+
\frac{1}{21}\\ \textrm{}\\
B_{8,F}(x)=x^{8}+21x^{7}+\frac{273}{2}x^{6}+364x^{5}+364x^{4}+\frac{273}{2}
x^{3}+21x^{2}+x+\frac{1}{36}\\ \textrm{}\\
B_{9,F}(x)=x^{9}+34x^{8}+357x^{7}+1547x^{6}+\frac{12376}{5}x^{5}+1547x^{4}+
357x^{3}+\\ \textrm{}\\+34x^{2}+x+\frac{1}{55}$
\end{enumerate}
\section{Fibonacci cobweb poset characterization}
\subsection{Fibonacci cobweb poset}

The Fibonacci cobweb poset $P$ has been invented by
A.K.Kwa\'sniewski in \cite{4.1,3.1,10.1} for the purpose of
finding combinatorial interpretation of fibonomial coefficients
and
eventually their reccurence relation. 

In \cite{4.1} A. K. Kwa\'sniewski defined cobweb poset $P$ as
infinite labeled digraph oriented upwards as follows: Let us label
vertices of $P$ by pairs of coordinates: $\langle i,j \rangle \in
{\bf N_{0}}\times {\bf N_{0}}$, where the second coordinate is the
number of level in which the element of $P$ lies (here it is the
$j$-th level) and the first one is the number of this element in
his level (from left to the right), here $i$. Following \cite{4.1}
we shall refer to $\Phi_{s}$ as to the set of vertices (elements)
of the $s$-th level, i.e.:
$$\Phi_{s}=\left\{\langle j,s \rangle ,\;\;1\leq j \leq F_{s}\right\},\;\;\;s\in{\bf N}\cup\{0\},$$
where $\{F_{n}\}_{n\geq 0}$ stands for Fibonacci sequence.

Then $P$ is a labeled graph $P=\left(V,E\right)$ where
$$V=\bigcup_{p\geq 0}\Phi_{p},\;\;\;E=\left\{\langle \,\langle j,p\rangle,\langle q,p+1
\rangle\,\rangle\right\},\;\;1\leq j\leq F_{p},\;\;1\leq q\leq
F_{p+1}.$$

 We can now define the partial order relation on $P$ as follows:
let\\ $x=\langle s,t\rangle, y=\langle u,v\rangle $ be elements of
cobweb poset $P$. Then
$$ ( x \leq_{P} y) \Longleftrightarrow
 [(t<v)\vee (t=v \wedge s=u)].$$


\subsection{DAG $\longrightarrow$ oDAG problem } In \cite{p.1} A. D.
Plotnikov considered the so called "DAG $\longrightarrow$ oDAG
problem". He determined condition when a digraph $G$ may  be
presented by the corresponding {\em dim } 2 poset $R$ and he
established the algorithm for finding it.

Before citing Plotnikov's results lat us recall  (following
\cite{p.1})  some indispensable definitions.

If $P$ and $Q$ are partial orders on the same set $A$, $Q$ is said
to be an {\bf extension} of $P$ if $a\leq_{P} b$ implies
$a\leq_{Q} b$, for all $a,b\in A$. A poset $L$ is a {\bf chain},
or a {\bf linear order} if we have either $a\leq_{L} b$ or
$b\leq_{L} a$ for any $a,b\in A$. If $Q$ is a linear order then it
is a {\bf linear extension} of $P$.

The {\bf dimension} $dim\ R$ of $R$ being a partial order is the
least positive integer $s$ for which there exists a family $F=(L_1
,L_2 ,\ldots,L_s)$ of linear extensions of $R$ such that $R=
\bigcap_{i=1}^{s} L_{i}$. A family $F=(L_1,L_2,\ldots,L_s)$ of
linear orders on $A$ is called a {\bf realizer} of $R$ on $A$ if
\[
R=\bigcap_{i=1}^{s} L_{i}.
\]

We denote by $D_{n}$ the set of all acyclic directed $n$-vertex
graphs without loops and multiple edges. Each digraph ${\vec
G}=(V,{\vec E})\in D_{n}$ will be called {\bf DAG}.

A digraph ${\vec G}\in D_{n}$ will be called {\bf orderable
(oDAG)} if there exists are $dim\ 2$ poset such that its Hasse
diagram coincide with the digraph ${\vec G}$.

Let  ${\vec G}\in D_{n}$ be a digraph, which does not contain the
arc $(v_{i},v_{j})$ if there exists the directed path
$p(v_{i},v_{j})$ from the vertex $v_{i}$ into the vertex $v_{j}$
for any $v_{i}$, $v_{j}\in V$. Such digraph is called {\bf
regular}. Let $D\subset D_{n}$ is the set of all regular graphs.

Let there is a some regular digraph ${\vec G}=(V,E)\in D$, and let
the chain ${\vec X}$ has three elements $x_{i_{1}}$, $x_{i_{2}}$,
$x_{i_{3}}\in X$ such that $i_{1}<i_{2}<i_{3}$, and, in the
digraph ${\vec G}$, there are not paths $p(v_{i_{1}},v_{i_{2}})$,
$p(v_{i_{2}},v_{i_{3}})$ and there exists a path
$p(v_{i_{1}},v_{i_{3}})$. Such representation of graph vertices by
elements of the chain ${\vec X}$ is called the representation in
{\bf inadmissible form}. Otherwise, the chain ${\vec X}$ presets
the graph vertices in {\bf admissible form}.

 Plotnikov showed that:

\begin{lemma}{\em \cite{p.1}}\label{l1}
\label{l1} A digraph ${\vec G}\in D_{n}$ may be represented by a
$dim\ 2$ poset if:
\renewcommand{\labelenumi}{(\arabic{enumi})}
\begin{enumerate}
\item there exist two chains ${\vec X}$ and ${\vec Y}$, each of
which is a linear extension of ${\vec G}_{t}$; \item the chain
${\vec Y}$ is a modification of ${\vec X}$ with inversions, which
remove the ordered pairs of ${\vec X}$ that there do not exist in
${\vec G}$.
\end{enumerate}
\end{lemma}
Above lemma results in the algorithm for finding {\em dim} 2
representation of a given DAG (i.e. corresponding oDAG) while the
following theorem establishes the conditions for constructing it.
\begin{thm}{\em \cite{p.1}}\label{t1}
\label{t1} A digraph ${\vec G}=(V,{\vec E})\in D_{n}$ can be
represented by $dim\ 2$ poset iff it is regular and its vertices
can be presented by the chain ${\vec X}$ in admissible form.
\end{thm}
\subsection{Fibonacci cobweb poset as DAG and oDAG} In this section
we show that Fibonacci cobweb poset is a DAG and it is orderable
(oDAG).

Obviously, cobweb poset $P=(V, E)$ defined above is a DAG (it is
directed acyclic graph without loops and multiple edges). One can
also verify that it is regular. For two elements $\langle i,
n\rangle , \langle j,m\rangle \in V$ a directed path $p(\langle i,
n\rangle , \langle j,m\rangle)\notin E$ will esist iff $n<m+1$ but
then  $(\langle i, n\rangle , \langle j,m\rangle)\notin E$ i.e.
$P$ does not contain the edge $(\langle i, n\rangle , \langle
j,m\rangle)$.

It is also possible to verify that vertices of cobweb poset $P$
can be presented in admissible form by the chain ${\vec X}$ being
a linear extension of cobweb $P$ as follows:
\begin{multline*}{\vec X}=\Big(\langle 1,0\rangle,\langle 1,1\rangle ,
\langle 1,2\rangle, \langle 1,3\rangle, \langle 2,3\rangle,
\langle 1,4\rangle, \langle 2,4\rangle, \langle 3,4\rangle,\langle
1,5\rangle, \langle 2,5\rangle, \langle 3,5\rangle,\\\langle
4,5\rangle,\langle 5,5\rangle,...\Big),\end{multline*}
 where

$$ ( \langle s,t\rangle \leq_{{\vec X}} \langle u,v\rangle) \Longleftrightarrow
 [(s\leq u)\wedge (t\leq v)]$$
 for $1\leq s \leq F_{t},\; 1\leq u \leq F_{v},\;\;\;t, v \in {\bf N}\cup\{0\}.$

 Fibonacci cobweb poset $P$ satisfies the conditions of Theorem
 \ref{t1} so it is oDAG. To find the chain ${\vec Y}$ being
a linear extension of cobweb $P$ one uses Lemma \ref{l1} and
arrives at: \begin{multline*} {\vec Y}=\Big(\langle
1,0\rangle,\langle 1,1\rangle , \langle 1,2\rangle, \langle
2,3\rangle, \langle 1,3\rangle, \langle 3,4\rangle, \langle
2,4\rangle, \langle 1,4\rangle,\langle 5,5\rangle, \langle
4,5\rangle, \langle 3,5\rangle,\\\langle 2,5\rangle,\langle
1,5\rangle,...\Big),\end{multline*}
 where

$$ ( \langle s,t\rangle \leq_{{\vec Y}} \langle u,v\rangle) \Longleftrightarrow
 [(t < v)\vee (t=v \wedge s\geq u)]$$
 for $1\leq s \leq F_{t},\; 1\leq u \leq F_{v},\;\;\;t, v \in {\bf
 N}\cup\{0\}$   and finally
 $$ (P,\leq_{P})={\vec X}\cap{\vec Y}.$$
\begin{rem}{\em For any sequence $\{a_{n}\}$ of
natural numbers one can define corresponding cobweb poset as
follows \cite{1c}:
$$\Phi_{s}=\left\{\langle j,s \rangle ,\;\;1\leq j \leq a_{s}\right\},\;\;\;s\in{\bf N}\cup\{0\},$$
and  $P=\left(V,E\right)$ where
$$V=\bigcup_{p\geq 0}\Phi_{p},\;\;\;E=\left\{\langle \,\langle j,p\rangle,\langle q,p+1
\rangle\,\rangle\right\},\;\;1\leq j\leq a_{p},\;\;1\leq q\leq
a_{p+1}$$
 with the partial order relation on $P$ :
$$ ( x \leq_{P} y) \Longleftrightarrow
 [(t<v)\vee (t=v \wedge s=u)]$$
 for  $x=\langle s,t\rangle, y=\langle u,v\rangle $ being elements of
cobweb poset $P$. Similary as above one can show that the family
of cobweb posets consist of DAGs representable  by corresponding
{\em dim } 2 posets (i.e. of oDAGs).}
\end{rem}

{\bf {Acknowledgements}}

I would like to thank Professor A. Krzysztof Kwa´sniewski for his
very helpful comments, suggestions, improvements and corrections
of this note.

\end{document}